\documentclass{amsart}
\title[Compactified stability manifold of a generic K3]{The Thurston compactification of the stability manifold of a generic analytic K3 surface}
\author{Anand Deopurkar}
\address{Mathematical Sciences Institute, Australian National University, Ngunnawal and Ngambri Country, Canberra, ACT, Australia}
\usepackage[margin=3cm]{geometry}
\usepackage{hyperref}
\usepackage{cleveref}
\usepackage{tikz}
\usepackage{tikz-cd}
\usepackage{amsthm}
\usepackage{thmtools}
\usepackage{cite}

\ifcsname theorem\endcsname{}\else\declaretheorem[parent=section]{theorem}\fi
\ifcsname corollary\endcsname{}\else\fi
\ifcsname lemma\endcsname{}\else\declaretheorem[sibling=theorem]{lemma}\fi
\ifcsname proposition\endcsname{}\else\declaretheorem[sibling=theorem]{proposition}\fi
\ifcsname conjecture\endcsname{}\else\fi
\ifcsname problem\endcsname{}\else\fi
\ifcsname question\endcsname{}\else\fi
\ifcsname definition\endcsname{}\else\fi
\ifcsname exercise\endcsname{}\else\fi
\ifcsname example\endcsname{}\else\fi
\ifcsname remark\endcsname{}\else\declaretheorem[sibling=theorem, style=remark]{remark}\fi

\providecommand {\Z}{{\bf Z}}

\providecommand{\Hom}{\operatorname{Hom}}

\providecommand{\rk}{\operatorname{rk}}

\begin{document}
\begin{abstract}
  Let \(X\) be an analytic K3 surface with \(\operatorname{Pic} X = 0\).
  We describe the closure of the Bridgeland stability manifold of \(X\) obtained using the masses of semi-rigid objects.
\end{abstract}

\maketitle
\section{Introduction}
Associated to a triangulated category \(\mathcal{C}\) is the complex manifold \(\operatorname{Stab}(\mathcal{C})\) whose points are the Bridgeland stability conditions on \(\mathcal{C}\) \cite{bri:07}.
Understanding the global geometry of \(\operatorname{Stab}(\mathcal{C})\) is an important question with far-reaching applications.
For example, when \(\mathcal{C}\) is the derived category of coherent sheaves on a K3 surface, the simple connectedness of \(\operatorname{Stab}(\mathcal{C})\) allows us to recover the group of auto-equivalences of \(\mathcal{C}\) \cite{bay.bri:17}.
When \(\mathcal{C}\) is the 2-Calabi--Yau category associated to a quiver, the topology of \(\operatorname{Stab}(\mathcal{C})\) has implications for the word/conjugacy problems and the \(K(\pi,1)\)-conjecture  for the associated Artin group \cite{qiu.woo:18,hen.lic:24}.

To better understand the global geometry of a non-compact space like \(\operatorname{Stab}(\mathcal{C})\), it is useful to have a compactification.
There have been several (partial) compactifications in the literature; see, for example, \cite{bro.pau.plo.ea:22,bar.mol.so:23,bol:23,bap.deo.lic:20}.
The goal of this paper is to completely describe the compactification constructed in \cite{bap.deo.lic:20} when \(\mathcal{C}\) is the derived category of coherent sheaves on a generic analytic K3 surface.

\begin{figure}[ht]
  \centering
  \begin{tikzpicture}[thick]
     \draw (120:1) -- (60:1) (120:1) -- (-90:1) (60:1) -- (-90:1)
      (120:1) -- (160:1) (160:1) -- (-90:1)
      (60:1) -- (20:1) (20:1) -- (-90:1)
      (160:1) -- (190:1) (190:1) -- (-90:1)
      (20:1) -- (-10:1) (-10:1) -- (-90:1)
      (-10:1) -- (-30:1)
      (190:1) -- (210:1);
      \draw (-100:1) circle (0.01) (-110:1) circle (0.01) (-120:1) circle (0.01);
      \draw (-80:1) circle (0.01) (-70:1) circle (0.01) (-60:1) circle (0.01);
      \draw [red, fill] (0,-1) circle (0.05);
    \end{tikzpicture}
    \caption{For an analytic K3 surface \(X\) with \(\operatorname{Pic}(X)= 0\), the compactified \(\mathbf{P} \operatorname{Stab}(X)\) is a closed disk, tiled by the translates of a triangle by the action of the spherical twist in \(\mathcal{O}_X\).  A distinguished point (red) in the boundary corresponds to the function \(\hom(\mathcal{O}_X, -)\).}
    \label{fig:disk}
  \end{figure}
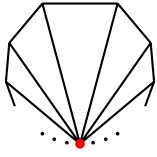

The compactification in \cite{bap.deo.lic:20} is motivated by viewing a stability condition as a metric, and in particular by Thurston's compactification of the Teichm\"uller space of hyperbolic metrics on a surface.
We recall the main idea.
Given a stability condition \(\sigma\) on \(\mathcal{C}\) and an object \(x \in \mathcal{C}\), the \emph{mass} of \(x\) with respect to \(\sigma\), denoted by \(m_{\sigma}(x)\), is the sum
\( m_{\sigma}(x) = \sum_i |Z_{\sigma}(x_i)|\),
where the \(x_i\) are the \(\sigma\)-Harder--Narasimhan (HN) factors of \(x\) and \(Z_{\sigma}\) is the central charge of \(\sigma\).
To construct the compactification, we fix a set of objects \(S\), and consider the map
\(m \colon \mathbf{P} \operatorname{Stab}(\mathcal{C}) = \operatorname{Stab}(\mathcal{C}) / \mathbf{C} \to \mathbf{P}^S\)
given by \(\sigma \mapsto [m_{\sigma}]\).
The proposed compactification is the closure of the image of \(m\).
\begin{theorem}\label{thm:main}
  Let \(X\) be an analytic K3 surface with \(\operatorname{Pic}(X) = 0\).
  Let \(S \subset D^b\operatorname{Coh} (X)\) be the set of semi-rigid objects.
  The map \(m \colon \mathbf{P} \operatorname{Stab}(D^b \operatorname{Coh}(X)) \to \mathbf{P}^S\) is a homeomorphism onto its image.
  The image is a 2-dimensional open ball and its closure is a 2-dimensional closed ball.
\end{theorem}
See \Cref{fig:disk} for an illustration of the compactified stability space.
The boundary contains a distinguished point represented by the function \(\hom(\mathcal{O}_X,-)\) (red point in \Cref{fig:disk}).
This is also the mass function of a lax stability condition in the sense of \cite{bro.pau.plo.ea:22}.
The other vertices are mass functions of lax pre-stability conditions.
The other boundary points do not have such interpretation.

  \Cref{thm:main} is a combination of \Cref{thm:homeo} and \Cref{prop:pi} in the main text.
The discussion of the points in the boundary is in \Cref{sec:boundary}.

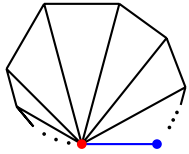
\begin{figure}[ht]
  \centering
  \begin{tikzpicture}[thick]
      \draw (120:1) -- (60:1) (120:1) -- (-90:1) (60:1) -- (-90:1)
      (120:1) -- (180:1) (180:1) -- (-90:1)
      (60:1) -- (20:1.2) (20:1.2) -- (-90:1)
      (180:1) -- (210:1) (210:1) -- (-90:1)
      (210:1) -- (230:1)      
      (20:1.2) -- (-10:1.4) (-10:1.4) -- (-90:1)
      (210:1) -- (230:1)
      (-10:1.4) -- (-20:1.4);
      \draw (-20:1.4) ++(-0.05,-0.1) circle (0.01) ++(-0.05,-0.1) circle (0.01) ++(-0.05,-0.1) circle (0.01);
      \draw (-100:1) circle (0.01) (-110:1) circle (0.01) (-120:1) circle (0.01);
      \draw[blue] (0,-1) -- (1,-1);
      \draw [red, fill] (0,-1) circle (0.05);
      \draw [blue, fill] (1,-1) circle (0.05);
    \end{tikzpicture}
      \caption{The closure of \(m_q(\mathbf{P} \operatorname{Stab}(X))\) is also a closed disk.
      The boundary has an additional interval, whose blue end-point is the \(q\)-hom functional \(\hom_q(\mathcal{O}_X,-)\).}
    \label{fig:q-disk}
  \end{figure}

  For a positive real number \(q\), the mass map has a natural \(q\)-analogue \(m_q\).
  The closure of the image of the stability manifold under \(m_q\) is also a closed disk.
  However, in its boundary, the red point in \Cref{fig:disk} is replaced by a closed interval (see \Cref{fig:q-disk}).

  For \(q = 1\), the distinguished point in the boundary has two interpretations: one as the hom function \(\hom(\mathcal{O}_X, -)\) and the second as the mass function of a lax stability condition \(\sigma\) in which \(\mathcal{O}_X\) is massless.
  For \(q \neq 1\), the two interpretations diverge.
  The \(q\)-hom function \(\hom_q(\mathcal{O}_X,-)\) yields the blue end-point in \Cref{fig:q-disk} and the \(q\)-mass function \(m_q(\sigma)\) yields the red end-point.
  
\begin{figure}[ht]
    \centering
    \begin{tikzpicture}[thick, scale=0.5]
      \draw[draw=none, fill=black, opacity=0.1] (-2,3) -- (-2,0) arc (180:0:1) arc (180:0:1) arc (180:0:1) (4,0) -- (4,3) -- (-2,3);
      \draw (-2,0) arc (180:0:1) (0,0) arc (180:0:1) (2,0) arc (180:0:1);
      \draw (-2,0) -- (-2,3) (0,0) -- (0,3) (2,0) -- (2,3) (4,0) -- (4,3);
      \draw (5,1) node {\(\cdots\)};
      \draw (-3,1) node {\(\cdots\)};
      \begin{scope}[xshift=12cm]
      \draw[draw=none, fill=black, opacity=0.1] (-2,3) -- (-2,0) arc (180:0:1) arc (180:0:0.8) arc (180:0:0.6) (2.8,0) -- (2.8,3) -- (-2,3);
      \draw (-2,0) arc (180:0:1) arc (180:0:0.8) arc (180:0:0.6);
      \draw (-2,0) -- (-2,3) (0,0) -- (0,3) (1.6,0) -- (1.6,3) (2.8,0) -- (2.8,3);
      \draw (3.4,1) node {\(\cdots\)};
      \draw (-3,1) node {\(\cdots\)};
      \draw[blue] (4,0) -- (4,3);      
      \end{scope}
    \end{tikzpicture}
    \caption{The tiling of the disk by triangles in the \(q = 1\) case (left) versus the \(q \neq 1\) case (right).}
    \label{fig:upper-half-plane}
  \end{figure}
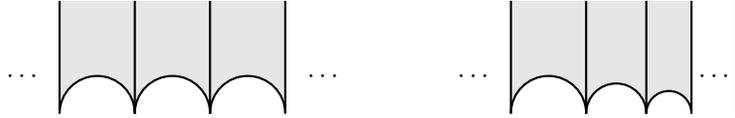

  We can reconcile the two pictures (\Cref{fig:disk} and \Cref{fig:q-disk}) by drawing them in the upper half plane instead of the disk (see \Cref{fig:upper-half-plane}).
  The \(q = 1\) picture (\Cref{fig:disk}) corresponds to the union of the translates of an ideal triangle by the transformation \(z \mapsto z + 1\).
  The only additional point in the closure (in the closed disk) is the point at infinity.
  The \(q \neq 1\) picture (\Cref{fig:q-disk}) corresponds to the union of the translates of an ideal triangle by the transformation \(z \mapsto q z + 1\).
  In this case, the closure (in the closed disk) contains an additional interval.
  This \(q\)-deformation is a simpler version of the \(q\)-deformed Farey tesselation observed in \cite{bap.bec.lic:22}.
  
In the course of the proof of the main theorem, we characterise all semi-rigid objects of \(D^b \operatorname{Coh} (X)\).
Up to twists by \(\mathcal{O}_X\) and homological shifts, the only such objects are the skyscraper sheaves \(\mathbf{k}_x\) (\Cref{prop:semirigids}).

There are a few other cases where the Thurston compactification of the stability manifold has been completely described.
These include: the 2-Calabi--Yau categories associated to quivers of rank 2 \cite{bap.deo.lic:20} and the derived categories of coherent sheaves on algebraic curves \cite{kik.kos.ouc:24}.
In \cite{kik.kos.ouc:25} the authors prove that for any (algebraic) K3 surface \(X\), taking \(S\) to be the set of spherical objects gives an injective map \(m \colon \mathbf{P} \operatorname{Stab}(X) \to \mathbf{P}^S\).
Understanding its image and its closure is an important goal.
The case of non-algebraic K3s treated here may be seen as a step towards it.

\subsection{Conventions}
An \emph{analytic K3 surface} is a connected, simply-connected, compact complex manifold \(X\) of dimension 2 with \(h^1(\mathcal{O}_X) = 0\).
By \(D^b(X)\) we mean the bounded derived category of the abelian category \(\operatorname{Coh}(X)\) of coherent sheaves on \(X\), as studied in \cite{huy.mac.ste:08}.
For a point \(x \in X\), we denote by \(\mathbf{k}_x\) the push-forward to \(X\) of the structure sheaf of \(x\), and call it the \emph{skyscraper sheaf} at \(x\).
By \(\operatorname{Stab}(X)\), we denote the set of (locally finite) Bridgeland stability conditions on \(D^b(X)\) with a numerical central charge; that is, where the central charge \(Z \colon K(D^b(X)) \to \mathbf{C}\) factors through the Chern character \(\operatorname{ch} \colon K(D^b(X)) \to H^{*}(X, \mathbf{Q})\).
We let \(\mathbf{P} \operatorname{Stab}(X)\) be the quotient of \(\operatorname{Stab}(X)\) by the standard action of \(\mathbf{C}\), in which \(z = x + i\pi y\) acts by scaling the central charge by \(e^{z}\) and shifting the slicing by \(y\).
Given a set \(S\), we let \(\mathbf{R}^S\) be the space of functions \(S \to \mathbf{R}\) with the product topology and \(\mathbf{P}^S\) the projective space \(\left(\mathbf{R}^S - \{0\}\right) / \text{scaling}\).
\subsection{Outline}
In \Cref{sec:stab}, we recall the description of stability conditions on an analytic K3 surface \(X\) with \(\operatorname{Pic} X = 0\).
In \Cref{sec:semirigid}, we identify the semi-rigid objects of \(D^b(X)\).
The bulk of the paper is \Cref{sec:mass}, in which we study the embedding of \(\mathbf{P} \operatorname{Stab}(X)\) given by the masses of semi-rigid objects.
In \Cref{sec:qmass}, we study the \(q\)-analogue of the mass embedding.
We do not include the definitions and the basic properties of stability conditions, and refer the reader to the original source \cite{bri:07} or exposition \cite{bay:11}.

\subsection{Acknowledgements}
This work is a part of a larger project with Asilata Bapat and Anthony Licata.
I am deeply grateful to have them as collaborators.
I thank Ziqi Liu, Emanuele Macri, Laura Pertusi, and Paolo Stellari for discussions related to this project.
I was supported by the Australian Research Council award \texttt{DP240101084}.

\section{Stability conditions on generic K3 surfaces}\label{sec:stab}
Throughout, fix an analytic K3 surface \(X\) with \(\operatorname{Pic} X = 0\).
Since \(X\) is a K3 surface, \(D^b(X)\) is a 2-Calabi--Yau  category.
That is, for \(x, y \in D^b(X)\), we have a natural isomorphism \[\Hom(x,y) \cong \Hom(y,x[2]).\]

\subsection{The Mukai lattice}
The Mukai lattice \(\mathcal{N}(X)\) of \(X\) is given by
\[ \mathcal{N}(X) = (H^0 \oplus H^4) (X, \mathbf{Z}).\]
Taking the class of \(X\) as a generator of the \(H^0\) summand and the class of a point \(x \in X\) as a generator of the \(H^4\) summand, we get an identification
\[ \mathcal{N}(X) = \mathbf{Z} \oplus \mathbf{Z}.\]
The Mukai pairing is then given by
\[ (\alpha_1,\alpha_2) \cdot (\beta_1, \beta_2) = \alpha_1 \beta_2 + \alpha_2\beta_1.\]
Given \(F \in D^b(X)\), we let \([F] = (\operatorname{ch}_0F, \operatorname{ch_0}F-\operatorname{ch}_2F) \in \mathcal{N}(X)\) be its Mukai vector.
Then we have
\[ [\mathcal{O}_X] = (1,1) \text{ and } [\mathbf{k}_x] = (0,1).\]
In particular, \([\mathcal{O}_X]\) and \([\mathbf{k}_x]\) form a basis of \(\mathcal{N}(X)\).
\subsection{Standard stability conditions}
We recall basic facts about stability conditions on \(X\) from \cite[\S~4]{huy.mac.ste:08}.
Let \(\mathcal{F}\) and \(\mathcal{T}\) be the full-subcategories of \(\operatorname{Coh}(X)\) consisting of torsion free and torsion sheaves, respectively.
Then \((\mathcal{F}, \mathcal{T})\) forms a torsion pair.
Let \(\mathcal{A}\) be the tilt of \(\operatorname{Coh}(X)\) in this torsion pair.
Explicitly,
\[  \mathcal{A} = \{E \in D^b(X) \mid H^{-1}(E) \in \mathcal{F} \text{ and } H^0(E) \in \mathcal{T} \text{ and for all }i \not \in \{0,1\}: H^i(E) = 0 \}.\]
Then \(\mathcal{A}\) is the heart of a bounded t-structure on \(D^b(X)\).

Let \(\mathbf{H} \subset \mathbf{C}\) be the (open) upper half plane.
As proved in \cite[\S~4.2]{huy.mac.ste:08}, for every \(z \in \mathbf{H} \cup \mathbf{R}_{<0}\), we have a stability condition \(\sigma_{z}\) on \(D^b(X)\) whose \((0,1]\) heart is \(\mathcal{A}\) and whose central charge is given by
\[ Z \colon [\mathbf{k}_x] \mapsto -1 \text{ and } Z \colon [\mathcal{O}_X] \mapsto -z.\]
For every \(w \in -\mathbf{H}\), we have a stability condition \(\sigma_{w}\) on \(D^b(X)\) whose \((0,1]\) heart is \(\operatorname{Coh}(X)\) and whose central charge is given by
\[ Z \colon [\mathbf{k}_x] \mapsto -1 \text{ and } Z \colon [\mathcal{O}_X] \mapsto -w.\]
See \Cref{fig:standardstability} for a sketch of the two central charges.
\begin{remark}
  The combined domain of the parameters \(z\) and \(w\) in \cite[\S~4.2]{huy.mac.ste:08} is \(\mathbf{C} - \mathbf{R}_{\geq -1}\).
  For us, it is \(\mathbf{C} - \mathbf{R}_{\geq 0}\).
  The difference is due to a slight change in parametrisation.
  The central charge of \(\sigma_z\) in \cite[\S~4.2]{huy.mac.ste:08} sends \(\mathbf{k}_x\) to \(-1\) (same as ours)  and \(\mathcal{O}_X\) to \(-z-1\) (we send it to \(-z\)).
  So our parametrisation and the parametrisation in \cite[\S~4.2]{huy.mac.ste:08} are related by \(z \mapsto z+1\).
\end{remark}
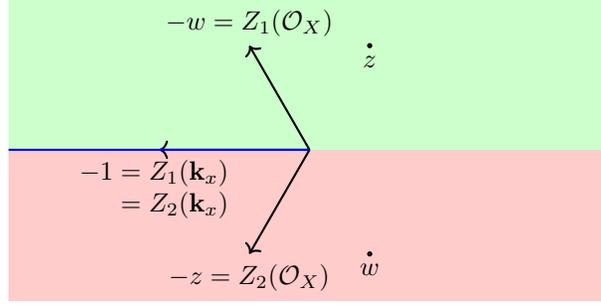
\begin{figure}[ht]
  \centering
  \begin{tikzpicture}[thick, scale=2]
    \draw[fill=red, opacity=0.2, draw=none] (2,0) -- (2,-1) -- (-2,-1) -- (-2,0) -- (2,0);
    \draw[fill=green, opacity=0.2, draw=none] (2,0) -- (2,1) -- (-2,1) -- (-2,0) -- (2,0);
    \draw[->] (0,0) --  (-1,0) node [below] {\parbox{6em}{\(-1 = Z_1(\mathbf{k}_x) \\ \phantom{-1} = Z_2(\mathbf{k}_x)\)}};
    \draw[->] (0,0) -- (-120:.8) node [below] {\(-z = Z_2(\mathcal{O}_X)\)};
    \draw[->] (0,0) -- (120:.8) node [above] {\(-w = Z_1(\mathcal{O}_X)\)};
    \draw[blue, opacity=1.0] (0,0) -- (-2,0);
    \draw (60:0.8) circle (0.01) node [below] {\(z\)};
    \draw (-60:0.8) circle (0.01) node [below] {\(w\)};
  \end{tikzpicture}
  \caption{For \(w \in -\mathbf{H}\) (red), a central charge \(Z_1\) as above defines a stability condition with heart \(\operatorname{Coh}(X)\). For \(z \in \mathbf{H}\) (green) and \(z \in \mathbf{R}_{<0}\) (blue), a central charge \(Z_2\) as above defines a stability condition whose heart is the tilt of \(\operatorname{Coh}(X)\) with respect to torsion and torsion-free sheaves.}
  \label{fig:standardstability}
\end{figure}

We call the stability conditions \(\sigma_z\) for \(z \in \mathbf{H} \cup -\mathbf{H} \cup \mathbf{R}_{<0}\) the \emph{standard stability conditions}.
We say that the stability conditions \(\sigma_z\) for \(z \in \mathbf{R}_{<0}\) are \emph{on the wall}, and the rest are \emph{off the wall}.

Let \(W_+\) (resp. \(W_-\) and \(W_0\)) be the union of the \(\mathbf{C}\)-orbits of the stability conditions \(\sigma_z\) for \(z \in \mathbf{H}\) (resp. \(-\mathbf{H}\) and \(\mathbf{R}_{<0}\)).
By definition, the sets \(W_+\), \(W_-\), and \(W_{0}\) are invariant under the \(\mathbf{C}\)-action.
It is easy to check that they are also invariant under the \(\widehat{\operatorname{GL}}_2^+(\mathbf{R})\)-action, and hence coincide with the sets with the same name defined in the proof of \cite[Theorem~4.8]{huy.mac.ste:08}.
Set \(W = W_+ \cup W_- \cup W_0\).

\subsection{All stability conditions}
Recall that the only spherical objects in \(D^b(X)\) are the shifts of \(\mathcal{O}_X\) (see \cite[Proposition~2.15]{huy.mac.ste:08}).
Let \(T \colon D^b(X) \to D^b(X)\) be the spherical twist in \(\mathcal{O}_X\).
\begin{proposition}\label{prop:Torbit}
  The set \(W \subset \operatorname{Stab}(X)\) is open and the union of its translates \(T^nW\), for \(n \in \mathbf{Z}\), is \(\operatorname{Stab}(X)\).
\end{proposition}
\begin{proof}
  That \(W\) is open is proved in \cite[Theorem~4.8]{huy.mac.ste:08}.
  That \(\operatorname{Stab}(X) = \bigcup T^nW\) is \cite[Corollary~4.7]{huy.mac.ste:08}.
\end{proof}

The following proposition allows us to identify the stability conditions in \(W_+\), \(W_{-}\), and \(W_0\).
Recall that since, up to shifts, \(\mathcal{O}_X\) is the only spherical object, it must be stable in any stability condition \cite[Proposition~2.15]{huy.mac.ste:08}.
\begin{proposition}\label{prop:type-comparison}
  Let \(\sigma\) be a stability condition and let \(\phi\) be the phase of \(\mathcal{O}_X\).
  Then \(\sigma\) is in \(W\) if and only if all the skyscraper sheaves \(\mathbf{k}_x\) are \(\sigma\)-stable of the same phase \(\psi\).
  In this case, we have
  \begin{enumerate}
  \item   \(\sigma \in W_{-}\) if \(\psi \in (\phi,\phi+1)\),
  \item  \(\sigma \in W_+\) if \(\psi \in (\phi+1,\phi+2)\),
  \item  \(\sigma \in W_0\) if \(\psi = \phi+1\).
  \end{enumerate}
\end{proposition}
\begin{proof}
  Since all skyscraper sheaves \(\mathbf{k}_x\) are \(\sigma\)-stable of the same phase for a standard stability condition, the same is true for any \(\sigma \in W\).
  Conversely, suppose all \(\mathbf{k}_x\) are \(\sigma\)-stable of the same phase.
  Using the \(\mathbf{C}\)-action, assume that their phase \(\psi\) is \(1\) and their central charge is \(-1\).
  By \cite[Proposition~4.6]{huy.mac.ste:08}, we conclude that \(\sigma\) is standard.

  Suppose \(\sigma = \sigma_z\) for \(z \in -\mathbf{H} \cup \mathbf{H} \cup \mathbf{R}_{<0}\).
  Whether \(z \in -\mathbf{H}\) or \(\mathbf{H}\) or \(\mathbf{R}_{<0}\) is distinguished by the phase \(\phi\) of \(\mathcal{O}_X\).
  For \(z \in -\mathbf{H}\), we have \(\phi \in (0,1)\);
  for \(z \in \mathbf{H}\), we have \(\phi \in (-1,0)\);
  and for \(z \in \mathbf{R}_{<0}\), we have \(\phi = 0\).
\end{proof}

\begin{proposition}\label{prop:T12}
  We have \(T W_+ = W_{-}\) and \(T^{-1} W_- = W_{+}\).
\end{proposition}
\begin{proof}
  We prove that for a standard \(\sigma \in W_-\), we have \(T^{-1}(\sigma) \in W_{+}\), and for a standard \(\sigma \in W_{+}\), we have \(T(\sigma) \in W_-\).
  Then the proposition follows.
  
  Take a standard \(\sigma \in W_-\)  and let us prove that \(T^{-1}(\sigma) \in W_{+}\).
  Let \(\phi \in (0,1)\) be the phase of \(\mathcal{O}_X\).
  It is easy to check that the ideal sheaves \(I_x\) of points \(x \in X\) are \(\sigma\)-stable of the same phase \(\psi \in (0,\phi)\).
  Let \(x \in X\) be any point.
  Since \(\Hom^{*}(\mathcal{O}_X, \mathbf{k}_x) = \mathbf{C}\), we have the exact triangle
  \[ \mathcal{O}_X \xrightarrow{\textrm{ev}} \mathbf{k}_x \to T \mathbf{k}_x \xrightarrow{+1}.\]
  Therefore, \( T \mathbf{k}_x = I_x[1]\).
  So \(T \mathbf{k}_x\) is \(\sigma\)-stable of phase \(\psi+1\).
  Therefore, \(T^{-1}I_x [1] = \mathbf{k}_x\) is \(T^{-1}(\sigma)\)-stable of phase \(\psi + 1 \in (1, \phi+1)\).
  On the other hand, \(T^{-1} \mathcal{O}_X = \mathcal{O}_X[1]\) is \(T^{-1}(\sigma)\)-stable of phase \(\phi\), so \(\mathcal{O}_X\) is \(T^{-1}(\sigma)\)-stable of phase \(\phi-1\).
  We now apply \Cref{prop:type-comparison}.

  Now take a standard \(\sigma \in W_+\)  and let us prove that \(T(\sigma) \in W_{-}\).
  Let \(\phi \in (-1,0)\)  be the phase of \(\mathcal{O}_{X}\).
  The objects \(T^{-1} \mathbf{k}_x\) are \(\sigma\)-stable of phase \(\psi \in (\phi+1,1)\)  (see \cite[Remark~4.3 (i)]{huy.mac.ste:08}).
  Therefore, the skyscraper sheaves \(\mathbf{k}_x\) are \(T(\sigma)\)-stable of phase \(\psi \in (\phi+1,1)\).
  Since \(\mathcal{O}_X\) is \(\sigma\)-stable of phase \(\phi\), it is \(T(\sigma)\)-stable of phase \(\phi+1\).
  We again apply \Cref{prop:type-comparison}.
\end{proof}

We now turn to the topology of the set of standard stability conditions and the stability conditions in \(W\).
Let \(H \subset \operatorname{Stab}(X)\) be the set of standard stability conditions.
Let \(R = \mathbf{C} \setminus \mathbf{R}_{\geq 0}\).
We have a map \(R \to H\) given by \(z \mapsto \sigma_z\).
We also have the projection map \(H \to \mathbf{P}W = W/ \mathbf{C}\).
\begin{proposition}\label{prop:chart}
 The maps \(R \to H\) and \(H \to \mathbf{P} W\) are homeomorphisms.
\end{proposition}
\begin{proof}
  By definition, the map \(R \to H\) is a bijection.
  By the proof of \cite[Theorem~4.8]{huy.mac.ste:08} (part (ii)), the map \(R \to H\) is continuous.
  Its inverse is given by \(\sigma \mapsto -Z_{\sigma}(\mathcal{O}_{X})\), which is also continuous.
  So \(R \to H\) is a homeomorphism.

  By \Cref{prop:type-comparison}, the map \(H \to \mathbf{P} W\) is surjective.
  Owing to the normalisation of the phase and mass of \(\mathbf{k}_x\), it is also injective.
  It remains to prove that the inverse is continuous.
  We know that \(W\) is an open subset of \(\operatorname{Stab}(X)\).
  It is also \(\mathbf{C}\)-invariant, so \(\mathbf{P} W\) is an open subset of \(\mathbf{P} \operatorname{Stab}(X)\).
  Thus, the map \(\mathbf{P} W \to \mathbf{P} \Hom(\mathcal{N}(X), \mathbf{C})\) is a local homeomorphism.
  We have the commutative diagram
  \[
    \begin{tikzcd}
      R\ar[equal]{d} & H \ar{l}[above]{\sim} & \mathbf{P} W\ar{l}\ar{d} \\
      R & & \mathbf{P} \Hom(\mathcal{N}(X), \mathbf{C})\ar{ll},
    \end{tikzcd}
  \]
  where the bottom map is given by \(Z \mapsto  Z(\mathcal{O}_X)/Z(\mathbf{k}_x)\).
  Since this map is continuous, it follows that \(\mathbf{P} W \to H\) is continuous.
\end{proof}

\section{Semi-rigid objects}\label{sec:semirigid}
Recall that an object \(F\) in \(D^b(X)\) is \emph{semi-rigid} if
\[ \hom^{i}(F,F) =
  \begin{cases}
     1 & \text{if } i = 0\\
    2 &\text{if } i = 1\\
    1 &\text{if } i = 2, \text{ and }\\
    0 &\text{otherwise}.
  \end{cases}
\]
For example, for \(x \in X\), the skyscraper sheaf \(F = \mathbf{k}_x\)  and the ideal sheaf \(F = I_x\) are semi-rigid.
We now characterise the semi-rigid objects of \(D^b(X)\).
Recall that \(T \colon D^b(X) \to D^b(X)\) is the spherical twist in \(\mathcal{O}_{X}\).
\begin{proposition}\label{prop:semirigids}
  Let \(X\) be a K3 surface with \(\operatorname{Pic} X = 0\).
  Let \(F \in D^b(X)\) be semi-rigid.
  Then there exists \(x \in X\) and integers \(m, n\)  such that \(F \cong T^n \mathbf{k}_x[m]\).
\end{proposition}
We split the proof in two lemmas.
\begin{lemma}\label{lem:semistable-semirigid}
  Fix a standard stability condition \(\sigma \in W_-\).
  Let \(F \in D^b(X)\) be semi-rigid and \(\sigma\)-semi-stable.
  Then there exists \(x \in X\) such that \(F\) or \(T^{-1}F\) is a shift of \(\mathbf{k}_x\).
\end{lemma}
\begin{proof}
  Since \(F\) is semi-rigid, \([F]\cdot[F] = 0\) in \(\mathcal{N}(X)\).
  So \([F]\) is an integer multiple of \((0,1)\) or \((1,0)\).

  Suppose \([F]\) is a multiple of \((0,1)\).
  Since \([\mathbf{k}_x] = (0,1)\), after applying a shift, we may assume that \(F\) is \(\sigma\)-semi-stable of the same phase as \(\mathbf{k}_x\), namely \(1\).
  It is easy to check that the abelian category of \(\sigma\)-semi-stable objects of phase \(1\) is \(\mathcal F\), the category of torsion sheaves on \(X\).
  It is a finite length category whose simple objects are the skyscraper sheaves \(\mathbf{k}_x\).
  So \(F\) is an iterated extension of skyscraper sheaves.
  Since \(\hom^1(F,F) = 2\), the Mukai lemma \cite[Lemma~2.7]{huy.mac.ste:08} implies that \(F\) must simply be a skyscraper sheaf.
  
  Suppose \([F]\) is a multiple of \((1,0)\).
  Then \([T^{-1}F]\) is a multiple of \((0,1)\) and \(T^{-1}F\) is semi-stable with respect to \(\tau = T^{-1}\sigma\).
  By \Cref{prop:T12}, we have \(\tau \in W_{+}\).
  By applying a rotation, assume that \(\tau\) is standard.
  Then, after applying a shift, we may assume that \(T^{-1}F\) is \(\tau\)-semi-stable of the same phase as \(\mathbf{k}_x\), namely \(1\).
  Again, it is easy to check that the abelian category of \(\tau\)-semi-stable objects of phase \(1\) is \(\mathcal{F}\).
  We now proceed as before.  
\end{proof}

Given a stability condition \(\sigma\), denote by \(\phi_{\sigma}^+\) and \(\phi_{\sigma}^{-}\) the highest and lowest phases of the factors in the \(\sigma\)-HN filtration.
If \(\sigma\) is clear from the context, we omit the subscript.
\begin{lemma}\label{lem:phasereduction}
  Fix a standard stability condition \(\sigma \in W_-\).
  Let \(F \in D^b(X)\) be a semi-rigid object.
  There exists a non-negative integer \(n\) such that \(T^nF\) is \(\sigma\)-semi-stable.
\end{lemma}
The proof relies on a result from \cite{bap.deo.lic:23} that allows us to reduce the phase spread by applying a spherical twist.
For the convenience of the reader, we recall the statement.
The statement applies to a triangulated category \(\mathcal{C}\), linear of finite type over a field, with a dg enhancement, equipped with a stability condition \(\tau\).
By \(\phi^{\pm}\) we denote the highest/lowest phases in the HN filtration and by \(T_{x}\) the spherical twist with respect to a spherical object \(x\).
\begin{theorem}[{\cite[Theorem~3.5]{bap.deo.lic:23}}]
  \label{thm:algorithm}
  Let \(x\) be a \(\tau\)-stable spherical object of \(\mathcal{C}\), and assume that it is the only \(\tau\)-stable object of its phase.
  Let \(y\) be an object of \(\mathcal{C}\) with no endomorphisms of negative degree.
  Suppose \(\phi(x) = \phi^+(y)\) (resp. \(\phi^-(y)\)) and let \(y' = T_x(y)\) (resp. \(T^{-1}_x(y)\)).
  If \(\phi^+(y) - \phi^-(y) > 0\), then
  \[ \phi^+(y') - \phi^-(y') < \phi^+(y) - \phi^-(y).\]
\end{theorem}
We now turn to the proof of \Cref{lem:phasereduction}.
\begin{proof}[Proof of \Cref{lem:phasereduction}]
  Since \(F\) is semi-rigid, all stable factors of \(F\) are either spherical or semi-rigid, and only one stable factor is semi-rigid \cite[Proposition~2.9]{huy.mac.ste:08}.
  The only spherical object, up to shift, is \(\mathcal{O}_X\).
  By \Cref{lem:semistable-semirigid}, the only semi-stable semi-rigid objects, up to shift, are \(\mathbf{k}_x\) and \(T \mathbf{k}_x\).
  In particular, the phases of the HN factors of \(F\) lie in the discrete subset of \(\mathbf{R}\) given by \[\left(\phi_{\sigma}(\mathcal{O}_X) + \mathbf{Z} \right) \cup \left(\phi_{\sigma}(\mathbf{k}_x) + \mathbf{Z} \right) \cup \left(\phi_{\sigma}(T^{1}\mathbf{k}_x) + \mathbf{Z}\right).\]
  Therefore, there exists a discrete \(\Phi \subset \mathbf{R}\) such that for every semi-rigid object \(F\), we have \[\phi^+(F) - \phi^-(F) \in \Phi.\]

  If \(F\) itself is semi-stable, we simply take \(n = 0\).
  Otherwise, for some \(i \in \mathbf{Z}\), a stable HN factor of \(F\) of highest or lowest phase must be \(\mathcal{O}_X[i]\).
  We apply \Cref{thm:algorithm} with \(y = F\) and \(x = \mathcal{O}_X[i]\).
  Then for \(F' = T F\) or \(F' = T^{-1} F\), we have
  \[ \phi^+(F') - \phi^-(F') < \phi^{+}(F) - \phi^{-1}(F).\]
  By repeated applications of \Cref{thm:algorithm} and using that \(\phi^+ - \phi^-\) lies in the discrete set \(\Phi \subset \mathbf{R}\), we conclude that there exists an integer \(n\) such that \(T^nF\) is semi-stable.
\end{proof}

Having proved the two lemmas, we are ready to prove \Cref{prop:semirigids}---the only semi-rigid objects of \(D^b(X)\), up to twisting by \(\mathcal{O}_X\) and shifting, are the skyscraper sheaves \(\mathbf{k}_x\).
\begin{proof}[Proof of \Cref{prop:semirigids}]
  Combine \Cref{lem:semistable-semirigid} and \Cref{lem:phasereduction}.
\end{proof}

\section{The mass embedding}\label{sec:mass}
Recall that \(X\) is an analytic K3 surface with \(\operatorname{Pic} X = 0\).
Let \(S\) be the set of isomorphism classes of semi-rigid objects of \(D^b(X)\).
In this section, we describe the mass embedding
\[ m \colon \mathbf{P} \operatorname{Stab}(X) \to \mathbf{P}^S\]
and the closure of its image.

\subsection{HN filtration of semi-rigid objects}
To understand the mass embedding, we must understand the HN filtrations of the objects of \(S\).
By \Cref{prop:semirigids}, the objects of \(S\), up to shift, are \(T^n \mathbf{k}_x\) for \(x \in X\) and \(n \in \mathbf{Z}\).
For points \(x,y \in X\), the behaviour of \(T^n \mathbf{k}_x\) and \(T^n \mathbf{k}_y\) is entirely analogous to each other.
So we lose nothing by fixing a particular point \(x \in X\) and taking
\[ S = \{T^{n} \mathbf{k}_x \mid n \in \mathbf{Z}\}.\]
We may then write the points of \(\mathbf{P}^S\) as homogeneous vectors \([x_n \mid n \in \mathbf{Z}] = [\cdots : x_{-1} :x_0 :x_1 : \cdots]\).
In these coordinates, the spherical twist \(T\) acts as a shift.

We first treat HN filtrations with respect to off the wall stability conditions.
\begin{proposition}\label{prop:HNII}
  Let \(\sigma \in W_{-}\).
  Then the \(\sigma\)-HN factors of \(F = T^n \mathbf{k}_x\), in decreasing order of phase, are as follows.
  \begin{enumerate}
  \item For \(n = 0\) and \(1\), the object \(F\) is stable.
  \item For \(n \geq 2\), the semi-stable (= stable) factors of \(F\) are \(T \mathbf{k}_x\) and \(\mathcal{O}_X[i]\) for \(0 \geq i \geq -n+2\).
  \item For \(n \leq -1\), the semi-stable (= stable) factors of \(F\) are \(\mathcal{O}_X[i]\) for \(-n \geq i \geq 1\) and \(\mathbf{k}_x\).
  \end{enumerate}
\end{proposition}
\begin{proof}
  Recall that \(\mathbf{k}_x\) and \(T \mathbf{k}_x = I_x[1]\) are stable for stability conditions in \(W_-\).
  So (1) follows.

  Consider the triangle
  \begin{equation}\label{eqn:twist1}
    \Hom^{*}(\mathcal{O}_X, T^{n-1} \mathbf{k}_x) \otimes \mathcal{O}_X \to T^{n-1} \mathbf{k}_x \to T^n \mathbf{k}_x \xrightarrow{+1}.
  \end{equation}
  We have
  \begin{align*}
    \Hom^{*}(\mathcal{O}_X, T^{n-1} \mathbf{k}_x) &= \Hom^{*}(T^{-n+1} \mathcal{O}_{X}, \mathbf{k}_x)\\
                                                  &= \Hom^{*}(\mathcal{O}_X[n-1], \mathbf{k}_x) \\
    &= \mathbf{C}[-n+1].
  \end{align*}
  By substituting in \eqref{eqn:twist1} and shifting, we get
  \begin{equation}\label{eqn:keytriangle}
    T^{n-1}\mathbf{k}_x \to T^n \mathbf{k}_x \to \mathcal{O}_{X} [-n+2] \xrightarrow{+1}.
  \end{equation}

  Let us assume \(n \geq 2\), and induct on \(n\).
  Assume we know that the HN factors of \(T^{n-1} \mathbf{k}_x\) (in decreasing order of phase) are \(T\mathbf{k}_{x}\) followed by \(\mathcal{O}_X[i]\) for \(0 \geq i \geq -n+3\).
  Concatenating the HN filtration of \(T^{n-1} \mathbf{k}_x\) and the map \(T^{n-1}\mathbf{k}_x \to T^n \mathbf{k}_x\), we obtain a filtration of \(T^n \mathbf{k}_x\) whose factors are \(T \mathbf{k}_x\) and \(\mathcal{O}_X[i]\) for \(0 \geq i \geq -n+2\).
  Since these factors are stable and appear in decreasing order of phase, this must be the HN filtration of \(T^n \mathbf{k}_x\).
  The induction step is complete.

  Now let us assume \(n \leq -1\), and induct on \(-n\).
  Consider the triangle
  \begin{equation}\label{eqn:keytriangle2}
    \mathcal{O}_{X} [-n] \to T^{n}\mathbf{k}_x \to T^{n+1} \mathbf{k}_x \xrightarrow{+1},
  \end{equation}
  obtained by replacing \(n\) by \(n+1\) in \eqref{eqn:keytriangle} and shifting.
  Assume we know that the HN factors of \(T^{n+1} \mathbf{k}_x\) (in decreasing order of phase) are \(\mathcal{O}_X[i]\) for \(-n-1 \geq i \geq 1\) and \(\mathbf{k}_x\).
  By augmenting the HN filtration of \(T^{n+1} \mathbf{k}_x\) by the map \(\mathcal{O}_X[-n] \to T^n \mathbf{k}_x\), we obtain a filtration of \(T^n \mathbf{k}_x\) whose factors are \(\mathcal{O}_X[i]\) for \(-n \geq i \geq 1\) and \(\mathbf{k}_x\).
  Since these factors are stable and appear in decreasing order of phase, this must be the HN filtration of \(T^n \mathbf{k}_x\).
   The induction step is complete.
 \end{proof}

For stability conditions on the wall, the HN filtration degenerates as expected.
\begin{proposition}\label{prop:HN0}
  Let \(\sigma \in W_{0}\).
  Then the \(\sigma\)-HN factors of \(F = T^n \mathbf{k}_x\), in decreasing order of phase, are as follows.
  \begin{enumerate}
  \item For \(n = -1\), \(0\) and \(1\), the object \(F\) is semi-stable.
  \item For \(n \geq 2\), the semi-stable factors of \(F\) are \(T\mathbf{k}_x\) and \(\mathcal{O}_X[i]\) for \(0 \geq i \geq -n+2\).
  \item For \(n \leq -2\), the semi-stable factors of \(F\) are \(\mathcal{O}_X[i]\) for \(-n \geq i \geq 2\) and \(T^{-1}\mathbf{k}_x\).
  \end{enumerate}
\end{proposition}
\begin{proof}
  The proof is analogous to the proof of \Cref{prop:HNII}.
\end{proof}

\subsection{The mass map}
\begin{figure}[ht]
  \centering
  \begin{tikzpicture}[thick]
    \draw[->] (0,0) -- node[below] {\(a\)} (-1,0) node [left] {\(Z(\mathbf{k}_x)\)};
    \draw[->] (0,0) -- node[right] {\(b\)} (1.3,1.5) node [right] {\(-Z(T \mathbf{k}_x)\)};
    \draw[->] (0,0) -- node[left] {\(c\)} (0.3,1.5) node [above] {\(Z(\mathcal{O}_X)\)};
    \draw[->, dashed] (0,0) -- (1,0);
    \draw (30:0.4) node {\tiny \(\theta\)};

    \begin{scope}[xshift=6cm]
       \draw[->] (0,0) -- node[below] {\(b\)} (-1,0) node [left] {\(Z(\mathbf{k}_x)\)};
       \draw[->] (0,0) -- node[right] {\(c\)} (1.3,1.5) node [right] {\(-Z(\mathcal{O}_X)\)};
       \draw[->] (0,0) -- node[left] {\(a\)} (0.3,1.5) node [above] {\(Z(T^{-1} \mathbf{k}_x)\)};
    \draw[->, dashed] (0,0) -- (1,0);
    \draw (30:0.4) node {\tiny \(\theta\)};
    \end{scope}
  \end{tikzpicture}
  \caption{We can use the cosine rule to reconstruct the central charge of a standard \(\sigma \in W_-\) from the masses \(a,b,c\) of \(\mathbf{k}_x, T\mathbf{k}_x, \mathcal{O}_X \) (left) and of \(\sigma \in W_+\) from the masses \(a,b,c\) of \(T^{-1} \mathbf{k}_x, \mathbf{k}_x, \mathcal{O}_X\) (right).}
  \label{fig:cosrule}
\end{figure}
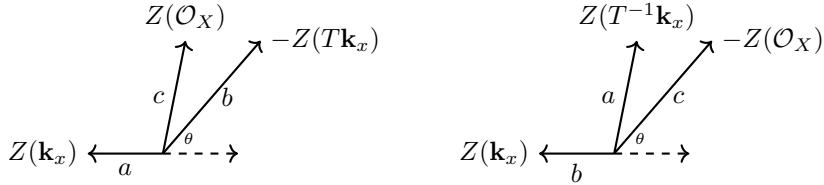

We now have the tools to describe the mass map
\[ m \colon \mathbf{P} \operatorname{Stab}(X) \to \mathbf{P}^S.\]
\begin{proposition}\label{prop:triangle0}
  Let \(\sigma \in \mathbf{P} W_{-}\).
  Set \(a = |Z_{\sigma}(\mathbf{k}_x)|\) and \(b = |Z_{\sigma}(T \mathbf{k}_{x})|\) and \(c = |Z_{\sigma}(\mathcal{O}_X)|\).
  \begin{enumerate}
  \item \label{eqn:triangleineqs}   The numbers \(a, b, c\) are positive real numbers satisfying
   \[
    b < a + c, \quad a < b+c, \quad c < a + b.
   \]
  \item We have
  \[ m_{\sigma} \colon T^n \mathbf{k}_x \mapsto
    \begin{cases}
      a - n c &\text{ if } n \leq 0,\\
      b + (n-1) c & \text{if } n \geq 1.
    \end{cases}
  \]
  \item Let \(\Delta_{0} \subset \mathbf{P}^S\) be the locally closed subset consisting of points of the form
   \[ [\cdots :a+2c:a+c:a:b:b+c:b+2c: \cdots],\]
   where \(a\) is at index 0 and \(b\) is at index \(1\), and where \(a,b,c\) are positive real numbers satisfying the inequalities in \eqref{eqn:triangleineqs}.
   Then \(m \colon \mathbf{P}W_{-} \to \Delta_{0}\) is a homeomorphism.
 \end{enumerate}
\end{proposition}
\begin{proof}
  Part (1) follows from the fact that the classes of \(\mathcal{O}_X\), \(\mathbf{k}_x\), and \(T \mathbf{k}_x\) satisfy
  \[ [\mathcal{O}_X] = [\mathbf{k}_x] - [T \mathbf{k}_x].\]

  Part (2) follows from \Cref{prop:HNII}.

  For part (3), note that \(\mathbf{P}W_- \to \Delta_0\) is continuous, so we must exhibit a continuous inverse.
  Let \(\Delta \subset \mathbf{P}^2\) be the set of points \([a:b:c]\) that satisfy the conditions in \eqref{eqn:triangleineqs}.
  Then we have a homeomorphism \(\Delta \to \Delta_0\) given by
  \[ [a:b:c] \mapsto [\cdots:a+2c:a+c:a:b:b+c:b+2c:\cdots].\]
  We use \([a:b:c] \in \Delta\) as coordinates on \(\Delta_0\).
  By \Cref{prop:chart}, the map \(w \mapsto \sigma_w\) gives a homeomorphism \(-\mathbf{H} \to \mathbf{P} W_{-}\).
  We use \(w \in -\mathbf{H}\) as a coordinate on \(\mathbf{P} W_-\).
  In these coordinates, writing down the inverse map \(\omega \colon \Delta_{0} \to \mathbf{P}W_-\) amounts to re-constructing the central charge given \(a,b,c\).
  This can be done using the cosine rule (see \Cref{fig:cosrule}).
  Precisely, we have
  \begin{equation}\label{eqn:inv-}
    \omega([a:b:c]) = - (b/a \exp(i \theta) - 1), \text{ where } \theta = \arccos \left( \frac{a^2+b^2-c^2}{2ab} \right) \in (0,\pi),
  \end{equation}
  which is the desired continuous inverse.
  \end{proof}

For \(n \in \mathbf{Z}\), let \(\Delta_n \subset \mathbf{P}^S\) be the locally closed subset consisting of points of the form
\[ [\cdots :a+2c:a+c:a:b:b+c:b+2c: \cdots],\]
where \(a\) is at index \(n\), and where \(a,b,c\) are positive real numbers satisfying the (strict) triangle inequalities.
Denote by \(T \colon \mathbf{P}^S \to \mathbf{P}^S\) the map that shifts the homogeneous coordinates rightwards by 1, so that \(\Delta_n = T^n \Delta_0\).
Recall that we also denote by \(T \colon \operatorname{Stab}(X) \to \operatorname{Stab}(X)\) the action of the spherical twist by \(\mathcal{O}_X\).
We have
\[ m (T (\sigma)) = T (m(\sigma)).\]
\Cref{prop:triangle0} implies that the mass map \(T^n \mathbf{P}W_- \to \Delta_n\) is a homeomorphism.
In particular, the mass map \(T^{-1} \mathbf{P}W_- = \mathbf{P}W_+ \to \Delta_{-1}\) is a homeomorphism.
It is useful to write the inverse \(\Delta_{-1} \to \mathbf{P}W_+\) using coordinates \([a:b:c]\) on \(\Delta_{-1}\) as in the proof of \Cref{prop:triangle0} and the coordinate on \(W_+\) given by \(z \in \mathbf{H}\).
Recall that the \([a:b:c]\) coordinates represent \(a = m(T^{-1} \mathbf{k}_x)\) and \(b = m(\mathbf{k}_x)\) and \(c = m(\operatorname{O}_X)\).
Then the map \([a:b:c] \mapsto z\) is (see \Cref{fig:cosrule}):
\begin{equation}\label{eqn:inv+}
  [a:b:c] \mapsto c/b \exp(i \theta), \text{ where } \theta = \arccos \left( \frac{b^2+c^2-a^2}{2bc} \right) \in (0,\pi).
\end{equation}

Let \(I_0 \subset \mathbf{P}^S\) be the set of points of the form
\[ [\cdots : a+2c :a+c: a : a+c : a+2c :\cdots], \]
where \(a\) is at index 0 and \(a,c\) are positive real numbers.
\begin{proposition}
  Let \(\sigma \in \mathbf{P}W_0\).
  Set \(a = |Z_{\sigma}(\mathbf{k}_x)|\) and \(c = |Z_{\sigma}(\mathcal{O}_X)|\).
  Then
  \[ m_{\sigma} \colon T^{n} \mathbf{k}_x \mapsto a + |n| c.\]
  Furthermore, the map \(m \colon \mathbf{P}W_0 \to I_0\) is a homeomorphism.
\end{proposition}
\begin{proof}
  The description of \(m_{\sigma}\) follows from \Cref{prop:HN0}.
  The inverse of \(m \colon \mathbf{P}W_0 \to I_0\) is given using the central charge \(Z(\mathbf{k}_x) = -1\) and \(Z(\mathcal{O}_X) = c/a\).
\end{proof}

\begin{proposition}\label{prop:domain}
  The map \(m \colon \mathbf{P}W \to \Delta_0 \cup I_0 \cup \Delta_{-1}\) is a homeomorphism.
\end{proposition}
See \Cref{fig:Wtriangle} for a sketch.
\begin{figure}
  \centering
    \begin{tikzpicture}[thick]
    \draw[fill=red, opacity=0.2, draw=none] (2,0) -- (2,-1) -- (-2,-1) -- (-2,0) -- (2,0);
    \draw[fill=green, opacity=0.2, draw=none] (2,0) -- (2,1) -- (-2,1) -- (-2,0) -- (2,0);
    \draw[blue, opacity=1.0] (0,0) -- (-2,0);
    \draw[white] (0,0) -- (2,0);

    \begin{scope}[xshift=4cm, xscale=0.8]
      \draw[fill=green, opacity=0.2, draw=none]  (0,0) -- (1,1) -- (2,0) -- (0,0);
      \draw[fill=red, opacity=0.2, draw=none]  (0,0) -- (1,-1) -- (2,0) -- (0,0);
      \draw[blue, opacity=1.0] (0,0) -- (2,0);
    \end{scope}

     \draw (2.5,0) edge [->] (3.5,0);
  \end{tikzpicture}
  \caption{The mass map gives a homeomorphism from the set of standard stability conditions parametrised by \(-\mathbf{H} \cup \mathbf{H} \cup \mathbf{R}_{<0}\) and the union of two open triangles and a segment that forms a wall between them.}
  \label{fig:Wtriangle}
\end{figure}
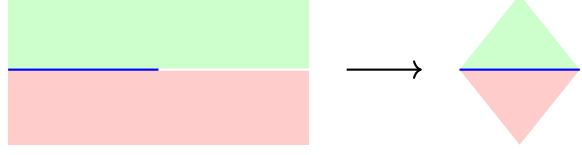
\begin{proof}
  The set \(\mathbf{P}W\) is the disjoint union of \(\mathbf{P}W_-\), \(\mathbf{P}W_{+}\), and \(\mathbf{P}W_0\).
  The sets \(\Delta_0\), \(I_{0}\), and \(\Delta_{-1}\) are also disjoint.
  Furthermore, the maps \(\mathbf{P}W_- \to \Delta_0\), \(\mathbf{P}W_+ \to \Delta_{-1}\), and \(\mathbf{P}W_0 \to I_0\) are homeomorphisms.
  So \(m \colon \mathbf{P} W \to \Delta_0 \cup I_0 \cup \Delta_{-1}\) is a continuous bijection.

  We check that the inverse is continuous.
  Since \(-\mathbf{H} \cup \mathbf{H} \cup \mathbf{R}_{<0} \to \mathbf{P}W\) is a homeomorphism, we use the former as local coordinates for \(\mathbf{P} W\).
  Let \(\overline \Delta \subset \mathbf{P}^2\) be the set of points \([a:b:c]\) where \(a,b,c\) are positive real numbers satisfying the triangle inequalities
  \[ b \leq a+c, \quad a < b+c, \quad c < a+b.\]
  It is easy to check that the map \(\overline \Delta \to \Delta_0 \cup I_0\) given by
 \[ [a:b:c] \mapsto [\cdots :a+c:a:b:b+c: \cdots]\]
 is a homeomorphism.
 So we may use \(a,b,c\) as local coordinates on \(\Delta_0 \cup I_0\).
  Using \eqref{eqn:inv-}, we see that the inverse map \(\Delta_0 \cup I_0 \to -\mathbf{H} \cup \mathbf{R}_{<0}\) is given in coordinates by
  \[ [a:b:c] \mapsto  - b/a \exp(i \theta) + 1, \text{ where } \theta = \arccos \left( \frac{a^2+b^2-c^2}{2ab} \right) \in [0,\pi),\]
  which is continuous.

  Let \(\overline \Delta' \subset \mathbf{P}^2\) be the set of points \([a:b:c]\) where \(a,b,c\) are positive real numbers satisfying the triangle inequalities
  \[ b < a+c, \quad a \leq b+c, \quad c < a+b.\]
  Then the map \(\overline \Delta' \to \Delta_{-1} \cup I_0\) given by
  \[ [a:b:c] \mapsto [\cdots :a+c:a:b:b+c: \cdots]\]
  is a homeomorphism.
  So we may use \(a,b,c\) as local coordinates on \(\Delta_{-1} \cup I_0\).
  Using \eqref{eqn:inv+}, we see that the inverse map \(\Delta_{-1} \cup I_0 \to \mathbf{H} \cup \mathbf{R}_{<0}\) is given in coordinates by
  \[
    [a:b:c] \mapsto c/b \exp(i \theta), \text{ where } \theta = \arccos \left( \frac{b^2+c^2-a^2}{2bc} \right) \in (0,\pi],
  \]
  which is continuous.

  Since the inverse is continuous on \(\Delta_0 \cup I_0\) and \(\Delta_{-1} \cup I_0\), we conclude that it is continuous on \(\Delta_0 \cup \Delta_{-1} \cup I_0\).
\end{proof}

Let \(D \subset \mathbf{P}^S\) be the union of the triangles \(\Delta_n\) for \(n \in \mathbf{Z}\) and the intervals \(I_n\) for \(n \in \mathbf{Z}\).
\begin{theorem}\label{thm:homeo}
  The mass map gives a homeomorphism \(m \colon \mathbf{P} \operatorname{Stab}(X) \to D\).
\end{theorem}
\begin{proof}
  By \Cref{prop:Torbit} and \Cref{prop:T12}, we see that \(\mathbf{P}\operatorname{Stab}(X)\) is the union of \(T^n \mathbf{P}W_-\) for \(n \in \mathbf{Z}\) and \(T^n \mathbf{P} W_0\) for \(n \in \Z\).
  From \Cref{prop:type-comparison}, it follows that this is a disjoint union.
  Likewise, \(D\) is the disjoint union of \(\Delta_n\) for \(n \in \mathbf{Z}\) and \(I_n\) for \(n \in \mathbf{Z}\).
  Since \(m \colon \mathbf{P}W_- \to \Delta_0\) and \(m \colon \mathbf{P}W_0 \to I_0\) are bijections, we conclude that \(m \colon \mathbf{P} \operatorname{Stab}(X) \to D\) is a bijection.
  It is also continuous.
  It remains to prove that the inverse is continuous.

  Let \(U = \Delta_0 \cup I_0 \cup \Delta_{-1}\).
  Observe that
  \[ U = \{[a_n] \in D \mid 2a_0 < a_1 + a_{-1}\}.\]
  So \(U \subset D\) is open.
  From \Cref{prop:domain}, we know that the inverse of \(m\) is continuous on \(U\).
  But \(T^n U\) for \(n \in \mathbf{Z}\) form an open cover of \(D\).
  So the inverse of \(m\) is continuous on \(D\).
\end{proof}

\subsection{Identifying the image and its closure}
Let \(\overline D \subset \mathbf{P}^S\) be the closure of \(D\).
Our next goal is to identify the homeomorphism classes of \(\overline D\) and \(D\).
To do so, it will be useful to work with an auxiliary space, which we now define.

Let \(I  \subset \mathbf{P}^{1}\) be the set of \([v:w]\) where \(v, w\) are non-negative real numbers.
Then \(I\) is homeomorphic to a closed interval.
Let \(\overline {\mathbf{R}} = \mathbf{R} \cup \{\pm \infty\}\) be the two point compactification of \(\mathbf{R}\), also homeomorphic to a closed interval.
Our auxiliary space will be \(\overline{\mathbf{R}} \times I\).

Define the transformation \(T\) on \(\overline{\mathbf{R}} \times I\) by
\[ T \colon (u,[v,w]) \mapsto (u+1,[v:w]).\]
Recall that we also denote by \(T\) the action of the spherical twist by \(\mathcal{O}_X\) on \(\mathbf{P} \operatorname{Stab}(X)\) and the rightward shift by \(1\) on \(\mathbf{R}^S\) and \(\mathbf{P}^S\).
(We intentionally use the same letter \(T\) to denote these maps, which are related.)
Our eventual goal is to understand \(\overline D\) via a \(T\)-equivariant parametrisation
\[ \pi \colon \overline{\mathbf{R}} \times I \to \overline D.\]
For \(n \in \mathbf{Z}\), set
\[ P_n = (\cdots , 2 , 1,  0 , 1 , 2 , \cdots) \in \mathbf{R}^S,\]
where the \(0\) is at index \(n\).
Note that \(P_n = T^nP_0\).
Set
\[ Q = (\cdots , 1, 1, 1, \cdots ) \in \mathbf{R}^S.\]
Observe that \(P_{0}, P_{1}\), and \(Q\) are the three vertices of the closure \(\overline{\Delta}_0\) of the triangle \(\Delta_0 \subset \mathbf{P}^S\), which is the homeomorphic image of \(\mathbf{P}W_-\).
The three sides of  \(\overline{\Delta}_0\) are the line segments \(P_{0}P_1\), \(P_{0}Q\), and \(P_1 Q\).
The open line segment \(P_0Q\) is the homeomorphic image of \(\mathbf{P}W_0\).
The entire picture is \(T\)-invariant, so the discussion above holds with \(0\) and \(1\) replaced by \(n\) and \(n+1\) for any \(n \in \mathbf{Z}\).

\begin{figure}[ht]
  \centering
  \begin{tikzpicture}[thick,yscale=0.5]
    \draw [red] (-3,1) -- (-3,-1);
    \foreach \i in {-2,...,2}{
    \draw (\i,1) -- (\i,-1);
    }
    \draw (-2.5,0) node {\(\cdots\)};
    \draw (2.5,0) node {\(\cdots\)};
    \draw (-3,1) -- (3,1);
    \draw [red] (-3,-1) -- (3,-1);
    \draw [red] (3,1) -- (3,-1);
    \begin{scope}[yshift=-5cm, yscale=2]
      \draw (120:1) -- (60:1) (120:1) -- (-90:1) (60:1) -- (-90:1)
      (120:1) -- (160:1) (160:1) -- (-90:1)
      (60:1) -- (20:1) (20:1) -- (-90:1)
      (160:1) -- (190:1) (190:1) -- (-90:1)
      (20:1) -- (-10:1) (-10:1) -- (-90:1)
      (-10:1) -- (-30:1)
      (190:1) -- (210:1);
      \draw (-100:1) circle (0.01) (-110:1) circle (0.01) (-120:1) circle (0.01);
      \draw (-80:1) circle (0.01) (-70:1) circle (0.01) (-60:1) circle (0.01);
      \draw [red, fill] (0,-1) circle (0.05);
    \end{scope}
    \draw (0,-1.5) edge [->] (0,-2.5);
    \draw (6,-3) node (S) {\(\overline{m(\mathbf{P} \operatorname{Stab}(X))} \subset \mathbf{P}^{S}\)};
    \draw (3.5,0) edge [->] (S);
    \draw (1.5,-5) edge [->] node[above] {\(\sim\)} (S);
  \end{tikzpicture}
  \caption{The map \(\pi \colon \overline{\mathbf{R}} \times [0,1] \to \mathbf{P}^S\) induces a homeomorphism from a closed disk \(B\) onto the closure of the image of \(\mathbf{P}\operatorname{Stab}(X)\).
    The disk \(B\) is obtained from the square \(\overline{\mathbf{R}} \times [0,1]\) by collapsing three sides (red).
    The \(\mathbf{Z}\)-indexed set of triangles tiling the image is the set of translates of the fundamental domain \(\overline{\mathbf{P}W_-}\) by the spherical twist \(T\).
  }
  \label{fig:accordion}
\end{figure}
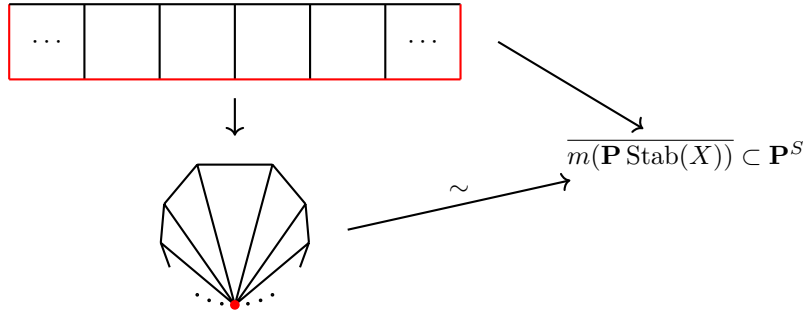

Consider the map \(\pi \colon [0,1] \times I \to \mathbf{P}^S\) defined by
\[ \pi(u,[v:w]) = (1-u)(wQ + vP_0) + u (wQ + vP_1).\]
Note that for  \(u = 0\) (resp. \(u = 1\)), the map \(\pi\) is a homeomorphism onto the closure of \(I_0\) (resp. \(I_1\)), which are the two sides \(P_{0}Q\) and \(P_1Q\) of the triangle \(\overline{\Delta}_0\).
For \(0 < u < 1\), the map \(\pi\) linearly interpolates between the two end-points \(\pi(0,[v:w])\) and \(\pi(1,[v:w])\), and hence its image is \(\overline \Delta_0\).
In fact, it is easy to check that the map
\[ \pi \colon [0,1] \times \left(I - \{[0:1]\}\right) \to \overline \Delta_0 - \{[\cdots 1 : 1: 1: \cdots]\}\]
is a homeomorphism, and \(\pi\) sends the entire segment \([0,1] \times [0:1]\) to the point \([\cdots 1 : 1: 1: \cdots]\).
Note that, with the \(T\) actions as before, we have
\[ T \pi(0,[v:w]) = \pi T(0,[v:w]).\]
Thus, \(\pi\) extends to a unique \(T\)-equivariant continuous map
\[ \pi \colon \mathbf{R} \times I \to \mathbf{P}^S.\]
Explicitly, for \(x = n + u\), where \(n \in \mathbf{Z}\) and \(u \in [0,1)\), we have
\[ \pi(x,[v:w]) = w Q + (1-u) v P_n + uvP_{n+1}.\]
Note that in \(\mathbf{R}^S\) we have
\[ \lim_{n \to \pm\infty} \frac{1}{n}P_n = (\cdots, 1, 1, 1, \cdots).\]
We extend \(\pi\) to \(\{\pm \infty\} \times I\) by setting
\[ \pi(\pm \infty, [v:w]) = [\cdots : 1 : 1 : 1 : \cdots].\]
Then, using the limit computation above, we see that \(\pi\) is continuous on \(\overline{\mathbf{R}} \times I\).
\begin{theorem}\label{prop:pi}
  The map \(\pi \colon \overline{\mathbf{R}} \times I \to \mathbf{P}^S\) is continuous.
  It sends the set \[C = \{\pm \infty\} \times I \cup \overline{\mathbf{R}} \times \{[0:1]\}\]
  to the point \([\cdots: 1: 1: 1: \cdots]\).
  Let \(\overline{\mathbf{R}} \times [0,1] \to B\) be the contraction of \(C\) to a point.
  Then the induced map \(\pi \colon B \to \mathbf{P}^S\) is a homeomorphism onto \(\overline D = \overline{m(\mathbf{P}\operatorname{Stab}(X))}\).
\end{theorem}
Note that \(B\) is homeomorphic to a closed disk.
See \Cref{fig:accordion} for a sketch.
\begin{proof}
  We have seen that \(\pi\) is continuous.
  It is easy to check that it is injective on the complement of \(C\), and its image on the complement of \(C\) does not include the point \([\cdots: 1:1:1: \cdots]\).
  It evidently sends all points of \(C\) to \([\cdots: 1:1:1: \cdots]\).
  So it induces a continuous injective map \(\pi \colon B \to \mathbf{P}^S\).
  Since \(B\) is compact and \(\mathbf{P}^S\) is Hausdorff, \(\pi\) maps \(B\) homeomorphically onto its image.
  By construction, \(\pi\) maps the interior of \(\overline{\mathbf{R}} \times I\) to \(D\).
  So \(\pi(B)\) must be \(\overline D\).
\end{proof}

\subsection{Points of the boundary}\label{sec:boundary}
Observe that \(\overline D\) contains the point \(\textcolor{red}{\bullet} = [\cdots :1:1:1: \cdots]\).
This is the common vertex (drawn in red in \Cref{fig:accordion}) of all the triangles that tile \(\overline D\).
It is the unique \(T\)-invariant point of \(\overline D\).
This point is precisely the projectivised hom function \(\hom(\mathcal{O}_X,-)\), whose value on \(T^n \mathbf{k}_x\) for any \(n \in \mathbf{Z}\) is
\[ \dim \Hom^{*}(\mathcal{O}_{X}, T^n\mathbf{k}_x) = 1.\]
The fact that \(\textcolor{red}{\bullet}\) is in the boundary is a reflection of the following more general fact.
\begin{theorem}[{\cite[Corollary~4.13]{bap.deo.lic:20}}]
  \label{thm:qhom}
  Let \(a\) be a spherical object of a triangulated category \(\mathcal{C}\), and assume that it is a stable object of a stability condition \(\sigma\).
  Let \(S\) be a set of objects of \(\mathcal{C}\) such that no object in \(S\) has an endomorphism of negative degree.
  For simplicity, also assume that no shift of \(a\) is in \(S\).
  Let \(T\) be the spherical twist in \(a\).
  Then, in \(\mathbf{P}^S\), we have the equality
  \[\lim_{n \to \pm \infty} T^n[m_{\sigma}] = [ \hom(a,-)].\]
\end{theorem}

The point \(\textcolor{red}{\bullet}\) also has an interpretation as the mass function of a lax stability condition in the sense of Broomhead, Pauksztello, Ploog, and Woolf \cite{bro.pau.plo.ea:22}.
We quickly recall the main features of the definition.
A \emph{lax pre-stability condition} consists of a slicing \(P\) and a compatible central charge \(Z\).
The central charge is allowed to vanish on the classes of non-zero semi-stable objects (such objects are called ``massless'').
A lax pre-stability condition is a \emph{stability condition} if it also satisfies the support property for massive stable objects.
That is, if there exists a choice of norm \(\|-\|\) on \(\mathcal{N}(X)\) and a positive constant \(c\) such that for every massive stable object \(s\), we have \(|Z(s)|/\|s\| \geq c\).
Given a lax pre-stability condition \(\sigma\) and an object \(x\), the mass of \(x\) with respect to \(\sigma\) is, as before, the sum \(\sum |Z_{\sigma}(x_i)|\) over the masses of HN pieces \(x_i\) of \(x\).

Recall that \(\mathcal{A}\) is the tilt of \(\operatorname{Coh} X\) in the torsion pair defined by torsion and torsion-free sheaves.
We let \(P\) to be the slicing defined by \(P(1) = \mathcal{A}\) and \(P(\phi) = 0\) for \(\phi \in (0,1)\).
The simple objects of \(P(1)\) are the skyscraper sheaves \(\mathbf{k}_x\) and the objects \(E[1]\), where \(E\) is a vector bundle on \(X\) with no non-trivial sub-bundles (see \cite[Remark~4.3 (iii)]{huy.mac.ste:08}).
We let \(Z(\mathcal{O}_X) = 0\) and \(Z(\mathbf{k}_x) = -1\).
\begin{proposition}\label{prop:red}
  The pair \((P,Z)\) as above defines a lax stability condition \(\sigma\) such that \[m(\sigma) = [\cdots : 1 : 1: 1 : \cdots] \in \mathbf{P}^S.\]
  Furthermore, \(\sigma\) is a limit of standard stability conditions.
\end{proposition}
\begin{proof}
  It is easy to check that the abelian category \(\mathcal{A}\) is of finite length (Noetherian and Artinian).
  So the slicing is locally finite.
  We now check the support property.
  Let \(E\) be a vector bundle with no non-trivial sub-bundles, and let \([E] = r[\mathcal{O}_X] + m[\mathbf{k}_x]\).
  Then \(r = \rk E\) and \(Z(E) = -m\).
  Assume that \(E\) is not isomorphic to \(\mathcal{O}_X\).
  Then \(\Hom(\mathcal{O}_X, E) = \Hom(E, \mathcal{O}_X) = 0\).
  So \[0 \geq \chi(\mathcal{O}_X, E) = 2r+m,\]
  and hence \(m \leq -2r\).
  With the norm on \(\mathcal{N}(X)\) in which \([\mathcal{O}_X]\) and \([\mathbf{k}_x]\) form an orthonormal basis, we see that \[|Z(E)|/\|E\| \geq \frac{|m|}{\sqrt{r^2+m^2}} \geq \frac{2}{\sqrt 5}.\]
  So the support property holds.

  The HN filtrations with respect to \(\sigma\) are the same as the HN filtrations with respect to a stability condition in \(\mathbf{P}W_0\).
  So the mass equality follows from \Cref{prop:HN0}.

  Finally, note that \(\sigma\) is the limit of stability conditions in \(W_0\) as \(Z(\mathcal{O}_{X})/Z(\mathbf{k}_x)\) approaches \(0\).
\end{proof}

Consider the points \(P_n\) of \(\overline D\).
These are the vertices of the tiling triangles other than the vertex \(\textcolor{red}{\bullet}\).
They form a single \(T\)-orbit, so it suffices to focus on one of them, say \(P_0 = [\cdots : 2:1:0:1:2: \cdots]\), with the \(0\) at index \(0\).
Note that this is the common vertex, other than \(\textcolor{red}{\bullet}\), of the triangles \(\mathbf{P}W_{+} \cong \Delta_{-1}\) and \(\mathbf{P}W_- = \Delta_0\).
This turns out to be the mass function of a lax pre-stability condition that does not satisfy the support property.
Let \(P\) be the same slicing as before, and set \(Z(\mathcal{O}_X) = 1\) and \(Z(\mathbf{k}_x) = 0\).
\begin{proposition}\label{prop:v0}
  The pair \((P,Z)\) as above defines a lax pre-stability condition \(\tau\) with
\[m(\tau) = [\cdots : 2: 1 : 0: 1 : 2: \cdots].\]
  The lax pre-stability condition \(\tau\) is a limit of standard stability conditions, and it does not satisfy the support property.
\end{proposition}
\begin{proof}
  Since \(\mathcal{A}\) is of finite length, the slicing is locally finite.
  So \(\tau\) is a lax pre-stability condition.
  The mass equality again follows from \Cref{prop:HN0}.
  Note that \(\tau\) is the limit of stability conditions in \(W_0\) as \(Z(\mathcal{O}_X)/Z(\mathbf{k}_x)\) approaches \(-\infty\).

  To see that the support property fails for \(\tau\), recall that the simple objects of \(\mathcal{A}\) are the skyscraper sheaves \(\mathbf{k}_x\) and shifts by 1 of vector bundles with no non-trivial sub-bundles.
  The skyscraper sheaves are massless, and hence do not obstruct the support property.
  However, for a fixed integer \(r \geq 2\) and sufficiently large \(m\) (depending on \(r\)), there exist vector bundles \(E\) of rank \(r\) on \(X\) with no non-trivial sub-bundles and \(c_2(E) = m\) (see \cite[Th\'eor\`eme~5.3]{ban.le-pot:87}).
  For such a vector bundle \(E\), we have \(Z_{\tau}(E) = r\) and \([E] = r [\mathcal{O}_X] + m [\mathbf{k}_x]\).
  By making \(m\) arbitrarily large, we can make \(|Z(E)|/\|E\|\) arbitrarily close to \(0\).
\end{proof}

Finally, consider a point on the open line segment joining \(P_{0}\) and \(P_{1}\).
This point is in the closure of \(\mathbf{P}W_- = \Delta_0\).
We claim that it is \emph{not} the mass function of a lax pre-stability condition arising as a limit of stability conditions \(W_-\).

To see this, it is helpful to consider a handful of other semi-stable objects.
Let \(n \geq m\) be positive integers.
Let \(x_1, \dots, x_n \in X\) be distinct points, and set \(R = \{x_1, \dots, x_n\}\).
We say that a morphism \( \pi \colon \mathcal{O}_X^{\oplus m} \to \mathcal{O}_R\) is \emph{generic} if for every subset \(T \subset R\), the induced map on global sections
\[ H^0(\mathcal{O}_X^{\oplus m}) \to H^0(\mathcal{O}_T)\]
has maximal rank, namely \(\min(m, |T|)\).

Let \(\sigma \in W_-\) be a standard stability condition.
Let \(I_{m,n}\) be the kernel of a generic morphism from \(\mathcal{O}_X^{\oplus m}\) to the structure sheaf of \(n\)-points.
Then it is easy to check that \(I_{m,n}\) is \(\sigma\)-stable.

Fix a point \(p \in \overline D\) on the line segment joining \(P_0\) and \(P_1\).
Then, for some \(t > 0\), we can write
\[ p = [ \cdots : 2 + t: 1: t : 1+2t : \cdots].\]
Take a sequence of standard stability conditions in \(W_{-}\) whose mass function approaches \(p\).
We claim that the sequence of slicings does not converge.

To see this, recall that the topology on the space of slicings is induced by the metric \(d\) defined as follows.
For a slicing \(P\) and non-zero object \(c\), let \(\phi^{\pm}_P(c)\) denote the highest/lowest phase of the \(P\)-HN factors of \(c\).
Then the distance \(d(P,Q)\) between two slicings \(P\) and \(Q\) is
\[ d(P,Q) = \operatorname{sup}_{c \neq 0}\left\{ \max(|\phi_{P}^+(c) - \phi_{Q}^+(c)|, |\phi^-_P(c) - \phi^-_Q(c)|)\right\}.\]
Suppose \(\tau\) is a lax stability condition that is a limit of a sequence of standard stability conditions in \(W_-\) with \(m(\tau) = p\).
Then \(\mathbf{k}_x\) and \(\mathcal{O}_X\) must be \(\tau\)-semi-stable.
Possibly after a rotation and scaling, we may assume that the central charge of \(\tau\) sends \(\mathbf{k}_x\) to \(-1\) and \(\mathcal{O}_X\) to \(-1-t\).
But then
\[ Z_{\tau}(I_{m,n}) = mZ_{\tau}(\mathcal{O}_X) - n Z_{\tau}(\mathbf{k}_x) = n - m(1+t).\]
It follows that for for every \((n,m)\) with \(n/m > (1+t)\), the sheaf \(I_{m,n}\) is \(\tau\)-semi-stable of phase \(0\) and for \(n/m < (1+t)\), it is \(\tau\)-semi-stable of phase \(1\).
But this is absurd.
Indeed, for a standard stability condition \(\sigma \in W_-\), we have
\[ \sup_{n/m > 1+t} \phi_{\sigma}(I_{n,m}) = \inf_{n/m < 1+t} \phi_{\sigma}(I_{n,m}),\]
so the same equality must hold in the limit.

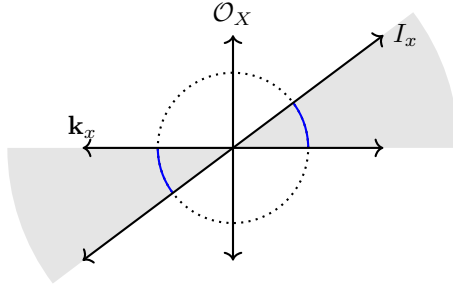
\begin{figure}
  \centering
  \begin{tikzpicture}[thick]
  \draw[draw=none, fill=gray, opacity=0.2] (0,0) -- (3,0) arc (0:37:3) (0,0);
  \draw[draw=none, fill=gray, opacity=0.2] (0,0) -- (-3,0) arc (180:217:3) (0,0);
  \draw[dotted] (0,0) circle (1);
  \draw[blue] (1,0) arc (0:37:1);
  \draw[blue] (-1,0) arc (180:217:1);
  \draw [->] (0,0) -- (0,1.5) node[above]{\(\mathcal{O}_{X}\)};
  \draw [->] (0,0) -- (0,-1.5);
  \draw [->] (0,0) -- (-2,0) node[above]{\(\mathbf{k}_x\)};
  \draw [->] (0,0) -- (2,0);
  \draw [->] (0,0) -- (2,1.5) node[right]{\(I_x\)};
  \draw [->] (0,0) -- (-2,-1.5);
\end{tikzpicture}
\caption{The central charges of semi-stable objects in a standard stability condition with heart \(\operatorname{Coh}X\) are the lattice points in the shaded region.  As a result, the phases are dense in the blue region of the unit circle.}
\label{fig:semistables}
\end{figure}

In summary, we see three kinds of behaviours on the boundary from the point of view of lax stability:
\begin{enumerate}
\item The object \(\mathcal{O}_X\) can become massless in a lax stability condition, leading to the limit mass function \(Q\).
\item The objects \(\mathbf{k}_x\) and \(I_x = T \mathbf{k}_x[-1]\) can become massless in a lax pre-stability condition, leading to the limit mass functions \(P_0\) and \(P_{1}\).
 \item Other semi-stable sheaves (for example, \(I_{m,n}\)) cannot become massless in lax pre-stability conditions.
\end{enumerate}
This trichotomy is consistent with the density of the phase diagram (discussed also in \cite[\S~12]{bro.pau.plo.ea:22}).
Let \(\sigma \in W_-\) be a standard stability condition.
It is easy to check that the classes \(r [\mathcal{O}_X] + n [\mathbf{k}_x]\) that support semi-stable sheaves are precisely
\(r = 0\) and \(n \geq 1\); or \(r \geq 1\) and \(n = 0\); or \(r \geq 1\) and \(-n \geq r\) (see \Cref{fig:semistables}).
Consider the phase diagram---the possible phases of semi-stable objects plotted on the unit circle.
There, \(\mathcal{O}_X\) is an isolated point, \(\mathbf{k}_x\) is a right accumulation point, and \(I_x\) is a left accumulation point.
At all points on the arc from \(\mathbf{k}_x[-1]\) to \(I_x\) (and its negative), the phase diagram is dense in the circle.
As \(\mathcal{O}_X\) becomes massless, the stability conditions converge preserving the support property.
As \(\mathbf{k}_x\) or \(I_x\) become massless, the slicings converge, but the support property is lost.
But if the central charge vanishes on a point on the open arc from \(\mathbf{k}_x[-1]\) to \(I_x\), even the slicings do not converge.

\section{The \(q\)-mass embedding}\label{sec:qmass}
Fix a positive real number \(q\).
Given a stability condition \(\sigma\) and an object \(x\), the \(q\)-mass of \(x\) with respect to \(\sigma\) is defined by
\[ m_{q, \sigma}(x) = \sum |Z_{\sigma}(x_i)| q^{\phi(x_i)},\]
where the sum is taken over the \(\sigma\)-HN factors \(x_i\) of \(x\), and \(\phi(x_i)\) is the phase of \(x_i\).
We have the map
\[ m_q \colon \mathbf{P} \operatorname{Stab}(X) \to \mathbf{P}^S\]
given by \(\sigma \mapsto m_{q,\sigma}\).
We describe the image of \(m_q\) and its closure for \(q \neq 1\).
Most of the arguments are direct analogues of the arguments for \(q = 1\), so we will be brief.

Let \(\sigma \in \mathbf{P}W_-\).
Set \(a = m_{q,\sigma}(\mathbf{k}_x)\) and \(b = m_{q,\sigma}(T\mathbf{k}_x)\) and \(c = m_{q,\sigma}(\mathcal{O}_X)\).
Owing to the triangle 
\[ \mathcal{O}_X \to \mathbf{k}_x \to T\mathbf{k}_x \xrightarrow{+1},\]
the positive real numbers \(a,b,c\) satisfy the \(q\)-triangle inequalities 
\begin{equation}\label{eqn:q-triangle-inequalities}
  b < a + qc, \quad a < b + c, \quad c < a + q^{-1}b.
\end{equation}
(See \cite[Proposition~3.3]{ike:21} for a proof of the \(q\)-triangle inequalities).
From the \(\sigma\)-HN filtration of \(T^n\mathbf{k}_x\) from \Cref{prop:HNII}, we get
\[ m_{q, \sigma} \colon T^n \mathbf{k}_x \mapsto
  \begin{cases}
    a + cq^{-n}+ \cdots + c q^{1} &\text{for \(n \leq -1\),}\\
    a &\text{for \(n = 0\),}\\
    b & \text{for \(n = 1\),}\\
    b + c q^0 + \cdots + c q^{-n+2} & \text{for \(n \geq 2\)}.
  \end{cases}
\]
So, in homogeneous coordinates, the \(q\)-mass map is
\[ m_q \colon \sigma \mapsto [ \cdots :a + cq+cq^2 : a + cq :a: b: b + c : b + c + c q^{-1} : \cdots ]\]

Let \(\Delta \subset \mathbf{P}^2\) be the set consisting of \([a:b:c]\) where \(a,b,c\) are positive real numbers satisfying \eqref{eqn:q-triangle-inequalities}.
Then the map \( \mathbf{P} W^- \to \Delta\)
that takes \(\sigma\) to \([m_{q,\sigma}(\mathbf{k}_x):m_{q,\sigma}(T\mathbf{k}_x): m_{q,\sigma}(\mathcal{O}_X)]\)
is a homeomorphism.
The proof is analogous to the proof of \Cref{prop:triangle0}~(3), but uses the \(q\)-analogue of the cosine rule \cite[Lemma~5.2]{bap.bec.lic:22}.

Consider \(\sigma \in \mathbf{P}W_0\).
With \(a,b,c\) as before, we have \(b = a + qc\).
From the \(\sigma\)-HN filtration of \(T^n\mathbf{k}_x\) from \Cref{prop:HN0}, we get
\[ m_{q, \sigma} \colon T^n \mathbf{k}_x \mapsto
  \begin{cases}
    a + cq^{-n}+ \cdots + c q^{1} &\text{for \(n \leq -1\),}\\
    a &\text{for \(n = 0\),}\\
    a + cq + \cdots + c q^{-n+2} & \text{for \(n \geq 1\)}.
  \end{cases}
\]
So, in homogeneous coordinates, the \(q\)-mass map is
\[ \sigma \mapsto [ \cdots :a + cq+cq^2 : a + cq :a: a+cq: a + cq + c: \cdots ].\]
Set \(I_{0} = m_q(\mathbf{P}W_0)\) and \(I_n = T^nI_0\).
Then \(m_q \colon T^n \mathbf{P}W_0 \to I_n\) is a homeomorphism.

Let \(D_q \in \mathbf{P}^S\) be the union of \(\Delta_n\) and \(I_n\) for \(n \in \mathbf{Z}\).
\begin{theorem}\label{thm:q-homeo}
  The \(q\)-mass map
  \[ m_q \colon \mathbf{P} \operatorname{Stab}(X) \to D_q\]
  is a homeomorphism.
\end{theorem}
The proof is analogous to the proof of \Cref{thm:homeo}.

We now identify the homeomorphism type of \(D_q\) and its closure \(\overline D_q\).
The basic technique is as before---by parametrising \(\overline D_q\) by a compactified infinite strip of squares.
But the resulting picture is slightly different.
Without loss of generality, assume \(q > 1\).

Recall that \(I  \subset \mathbf{P^1}\) is the set of \([v:w]\) where \(v, w\) are non-negative real numbers.
Let \(\overline {\mathbf{R}} = \mathbf{R} \cup \{\pm \infty\}\) be the two point compactification of \(\mathbf{R}\).
Define the transformation \(T\) on \(\overline{\mathbf{R}} \times I\) by
\[ T \colon (u,[v,w]) \mapsto (u+1,[qv:w]).\]
Recall that we also denote by \(T\) the action of the spherical twist by \(\mathcal{O}_X\) on \(\mathbf{P} \operatorname{Stab}(X)\) and the rightward shift by \(1\) on \(\mathbf{R}^S\) and \(\mathbf{P}^S\).

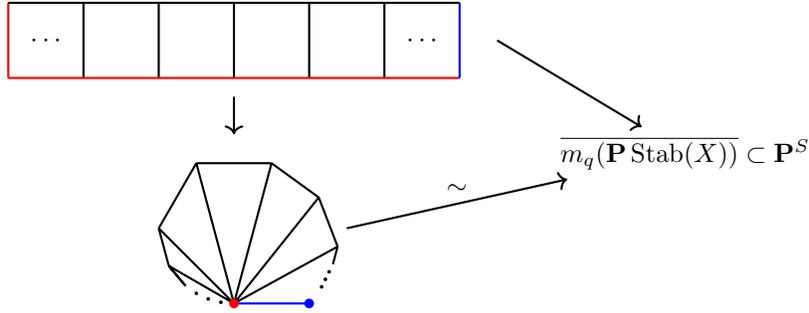
\begin{figure}[ht]
  \centering
  \begin{tikzpicture}[thick,yscale=0.5]
    \draw [red] (-3,1) -- (-3,-1);
    \foreach \i in {-2,...,2}{
    \draw (\i,1) -- (\i,-1);
    }
    \draw (-2.5,0) node {\(\cdots\)};
    \draw (2.5,0) node {\(\cdots\)};
    \draw (-3,1) -- (3,1);
    \draw [red] (-3,-1) -- (3,-1);
    \draw [blue] (3,1) -- (3,-1);
    \begin{scope}[yshift=-5cm, yscale=2]
      \draw (120:1) -- (60:1) (120:1) -- (-90:1) (60:1) -- (-90:1)
      (120:1) -- (180:1) (180:1) -- (-90:1)
      (60:1) -- (20:1.2) (20:1.2) -- (-90:1)
      (180:1) -- (210:1) (210:1) -- (-90:1)
      (210:1) -- (230:1)      
      (20:1.2) -- (-10:1.4) (-10:1.4) -- (-90:1)
      (210:1) -- (230:1)
      (-10:1.4) -- (-20:1.4);
      \draw (-20:1.4) ++(-0.05,-0.1) circle (0.01) ++(-0.05,-0.1) circle (0.01) ++(-0.05,-0.1) circle (0.01);
      \draw (-100:1) circle (0.01) (-110:1) circle (0.01) (-120:1) circle (0.01);
      \draw[blue] (0,-1) -- (1,-1);
      \draw [red, fill] (0,-1) circle (0.05);
      \draw [blue, fill] (1,-1) circle (0.05);
    \end{scope}
    \draw (0,-1.5) edge [->] (0,-2.5);
    \draw (6,-3) node (S) {\(\overline{m_{q}(\mathbf{P} \operatorname{Stab}(X))} \subset \mathbf{P}^{S}\)};
    \draw (3.5,0) edge [->] (S);
    \draw (1.5,-5) edge [->] node[above] {\(\sim\)} (S);
  \end{tikzpicture}
  \caption{The map \(\pi_{q} \colon \overline{\mathbf{R}} \times [0,1] \to \mathbf{P}^S\) induces a homeomorphism from a closed disk \(B\) onto the closure of the image of \(\mathbf{P}\operatorname{Stab}(X)\) under the \(q\)-mass map.
    The disk \(B\) is obtained from the square \(\overline{\mathbf{R}} \times [0,1]\) by collapsing two sides (red).
  }
  \label{fig:q-accordion}
\end{figure}

We define a \(T\)-equivariant parametrisation
\[ \pi_q \colon \overline{\mathbf{R}} \times I \to \overline D_q,\]
which is a \(q\)-analogue of the parametrisation \(\pi\) from \Cref{prop:pi}.
For \(n \in \mathbf{Z}\), set
\[ P_n = (\cdots , 1+q , 1,  0 , 1 , 1+q^{-1} , \cdots) \in \mathbf{R}^S,\]
where the \(0\) is at index \(n\).
Note that \(P_n = T^nP_0\).
Set
\[ Q = (\cdots , 1, 1, 1, \cdots ) \in \mathbf{R}^S.\]
Consider the map \(\pi_{q} \colon [0,1] \times I \to \mathbf{P}^S\) defined by
\[ \pi_q(u,[v:w]) = (1-u)(wQ + vP_0) + u (wQ + q^{-1}vP_1).\]
Note that for \(u = 0\) (resp. \(u = 1\)), the map \(\pi_q\) is a homeomorphism onto the closure of \(I_0\) (resp. \(I_1\)), which are the two sides of the triangle \(\overline{\Delta}_0\).
For \(0 < u < 1\), the map \(\pi_q\) linearly interpolates between the two end-points \(\pi_q(0,[v:w])\) and \(\pi_q(1,[v:w])\), and hence its image is \(\overline \Delta_0\).
Also observe that \(\pi_q(u,[0:w]) = Q\).
Furthermore, with the \(T\) actions as before, we have
\[ T \pi_q(0,[v:w]) = \pi_q T(0,[v:w]).\]
Thus, \(\pi_q\) extends to a unique \(T\)-equivariant continuous map
\[ \pi_q \colon \mathbf{R} \times I \to \mathbf{P}^S.\]
Explicitly, for \(x = n + u\), where \(n \in \mathbf{Z}\) and \(u \in [0,1)\), we have
\[ \pi_q(x,[v:w]) = w Q + (1-u) q^{-n}v P_n + u q^{-n-1}vP_{n+1}.\]
Let \(\delta = 1 + q^{-1} + q^{-2} + \cdots\).
Then, in \(\mathbf{R}^S\) we have
\[ \lim_{n \to -\infty} P_n = (\dots, \delta, \delta, \delta, \dots) \text{ and } \lim_{n \to \infty}q^{-n}P_n = (\dots, q \delta,\delta ,q^{-1} \delta,\dots).\]
(On the right hand side of the last equation, the \(\delta\) is at index \(-1\).)
Extend \(\overline \pi_q\) to \(\{\pm \infty\} \times I\) by setting
\[ \pi_q(-\infty, [v:w]) = [\cdots : 1 : 1 : 1 : \cdots],\]
and
\[ \pi_q(+\infty, [v:w]) = w [\cdots : 1 : 1 : 1 : \cdots] + v [\cdots : q : 1 : q^{-1} : \cdots].\]
Using the limit computation above, it follows that this extension is continuous on \(\overline{\mathbf{R}} \times I\).
\begin{theorem}\label{prop:q-pi}
  The map \(\pi_{q} \colon \overline{\mathbf{R}} \times [0,1] \to \mathbf{P}^S\) is continuous.
  It sends the set \[C = \{- \infty\} \times I \cup \overline{\mathbf{R}} \times \{[0:1]\}\]
  to the point \([\cdots: 1: 1: 1: \cdots]\).
  Let \(\overline{\mathbf{R}} \times [0,1] \to B\) be the contraction of \(C\) to a point.
  Then the induced map \(\pi_{q} \colon B \to \mathbf{P}^S\) is a homeomorphism onto \(\overline D_q = \overline{m_{q}(\mathbf{P}\operatorname{Stab}(X))}\).
\end{theorem}
The proof is analogous to that of \Cref{prop:pi}.
See \Cref{fig:q-accordion} for a sketch.

Instead of a unique \(T\)-fixed point of \(\overline D_q\), as was the case for \(q = 1\), for \(q \neq 1\) we have two such points.
These are the blue and red end-points of the blue interval in \Cref{fig:q-accordion}.
The blue end-point is the point \(\textcolor{blue}{\bullet} = [\cdots : q : 1 : q^{-1} : \cdots  ]\).
It is the \(q\)-hom function \(\hom_q(\mathcal{O}_X, -)\), whose value on \(T^n\mathbf{k}_x\) is
\[ \dim_q \Hom^{*}(\mathcal{O}_X, T^n \mathbf{k}_x) = q^{-n}.\]
(By definition, \(\dim_q\) of the graded vector space \(\mathbf{C}[m]\) is \(q^m\)).
Again, the fact that \(\hom_q(\mathcal{O}_X,-)\) appears in the closure of the \(q\)-mass embedding is a reflection of a general theorem---the \(q\)-analogue of \Cref{thm:qhom} (see \cite[Corollary~4.13]{bap.deo.lic:20}).

Note that \(\textcolor{blue}{\bullet}\) is not in the closure of the standard stability conditions \(\mathbf{P}W\), nor is it in the closure of \(T^n \mathbf{P}W\) for any fixed \(n\).
To reach \(\textcolor{blue}{\bullet}\), we must traverse an infinite sequence of hearts.
It is not the \(q\)-mass function of a lax pre-stability condition.

The red end-point is the point \(\textcolor{red}{\bullet} = [\cdots : 1 : 1 : 1 : \cdots]\).
It is the \(q\)-mass function of the lax stability condition \(\sigma\) from \Cref{prop:red}.

The other vertices of the triangles form one orbit, and are \(q\)-mass functions of lax pre-stability conditions.
For example, the vertex \(P_0 = [ \cdots: 1+q : 1 : 0 : 1 : 1+q^{-1} : \cdots ]\) is the \(q\)-mass function of the lax pre-stability condition \(\tau\) from \Cref{prop:v0}.

\bibliographystyle{amsalpha}
\providecommand{\bysame}{\leavevmode\hbox to3em{\hrulefill}\thinspace}
\providecommand{\MR}{\relax\ifhmode\unskip\space\fi MR }
\providecommand{\MRhref}[2]{\href{http://www.ams.org/mathscinet-getitem?mr=#1}{#2}
}
\providecommand{\href}[2]{#2}

\end{document}